\documentclass[final]{siamltex}

\usepackage{amssymb}
\usepackage{amsmath}
\usepackage{graphicx}
\usepackage{subfig}
\usepackage{algorithm}
\usepackage{algpseudocode}
\usepackage{psfrag}
\usepackage{url}
\usepackage{array}
\usepackage[sort]{natbib}

\usepackage[usenames]{color}
\usepackage[normalem]{ulem}

\newcolumntype{S}{>{\centering\arraybackslash} m{0.8cm}}
\newcolumntype{F}{>{\centering\arraybackslash} m{.2\linewidth}}

\begin{document}

\title{A Fourier-based Approach for Iterative 3D Reconstruction from Cryo-EM Images}
\author{Lanhui Wang\thanks{The Program in Applied and Computational Mathematics (PACM), Princeton University, Fine Hall, Washington Road, Princeton, NJ 08544-1000
Princeton University, \texttt{lanhuiw@math.princeton.edu}, Corresponding author. Tel.: +1 609 258 5785; fax: +1 609 258 1735. }  \and Yoel Shkolnisky\thanks{Department of Applied Mathematics, School of Mathematical Sciences, Tel Aviv University, Tel Aviv 69978, Israel, \texttt{yoelsh@post.tau.ac.il} }
\and Amit Singer\thanks{Department of Mathematics and PACM, Princeton University, Fine Hall, Washington Road, Princeton, NJ 08544-1000
Princeton University, \texttt{amits@math.princeton.edu} } }

\date{}
\maketitle

\begin{abstract}
A major challenge in single particle reconstruction methods using cryo-electron microscopy is to attain a resolution sufficient to interpret fine details in three-dimensional (3D) macromolecular structures. Obtaining high resolution 3D reconstructions is difficult due to unknown orientations and positions of the imaged particles, possible incomplete coverage of the viewing directions, high level of noise in the projection images, and limiting effects of the contrast transfer function of the electron microscope. In this paper, we focus on the 3D reconstruction problem from projection images assuming an existing estimate for their orientations and positions. We propose a fast and accurate Fourier-based Iterative Reconstruction Method (FIRM) that exploits the Toeplitz structure of the operator ${\bf A}^{*}{\bf A}$, where $\bf A$ is the forward projector and ${\bf A}^{*}$ is the back projector. The operator ${\bf A}^{*}{\bf A}$ is equivalent to a convolution with a kernel. The kernel is pre-computed using the non-uniform Fast Fourier Transform  and is efficiently applied in each iteration step. The iterations by FIRM are therefore considerably faster than those of traditional iterative algebraic approaches, while maintaining the same accuracy even when the viewing directions are unevenly distributed. The time complexity of FIRM is comparable to the direct Fourier inversion method. Moreover, FIRM combines images from different defocus groups simultaneously and can handle a wide range of regularization terms. We provide experimental results on simulated data that demonstrate the speed and accuracy of FIRM in comparison with current methods.
\end{abstract}


\section{Introduction}
Single particle reconstruction (SPR) from cryo-electron microscopy (cryo-EM) \cite{Frank1996,vanheel_cryo-em} is an emerging technique for determining the 3D structure of macromolecules.  One of the main challenges in SPR is to attain a resolution of 4\AA{}  or better, thereby allowing interpretation of atomic coordinates of macromolecular maps \cite{Frank_resolution, Fred_resolution}. Although X-ray crystallography and NMR spectroscopy can achieve higher resolution levels ($\sim1$\AA{} by X-ray crystallography and 2-5\AA{} by NMR spectroscopy), these traditional methods are often limited to relatively small molecules. In contrast, cryo-EM is typically applied to large molecules or assemblies with size ranging from 10 to 150 nm, such as ribosomes \cite{ribo_cryo-em}, protein complexes, and viruses.

Cryo-EM is used to acquire 2D projection images of thousands of individual, identical frozen-hydrated macromolecules at random unknown orientations and positions. The collected images are extremely noisy due to the limited electron dose used for imaging
to avoid excessive beam damage. In addition, the unknown pose parameters (orientations and positions) of the imaged particles need to be estimated for 3D reconstruction. An ab-initio estimation of the pose parameters using the random-conical tilt technique \cite{GWBP_radermacher} or common-lines based approaches \cite{VanHeel1987111, Amit_voting, Amit_eig_sdp} are often applied after multivariate statistical data compression \cite{Multivariatestat,vanHeel1981187} and classification techniques \cite{vanHeel1984165, GWBP_Penczek1992, Amit_classavg} that are used to sort and partition the large set of images by their viewing directions, producing ``class averages'' of enhanced signal-to-noise ratio (SNR). Using the ab-initio estimation of the pose parameters, a preliminary 3D map is reconstructed from the images by a 3D reconstruction algorithm. The initial model is then iteratively refined \cite{Refinement} in order to obtain a higher-resolution 3D reconstruction. In each iteration of the refinement process, the current 3D model is projected at several pre-chosen viewing directions and the resulting images are matched with the particle images, giving rise to new estimates of their pose parameters. The new pose parameters are then used to produce a refined 3D model using a 3D reconstruction algorithm. This process is repeated for several iterations until convergence. Clearly, a fast and accurate 3D reconstruction algorithm is needed for both the initial model reconstruction and for the refinement process. The focus of this paper is the 3D reconstruction problem with given pose parameters.

\begin{figure}
\begin{centering}
\includegraphics[width=0.53\paperwidth]{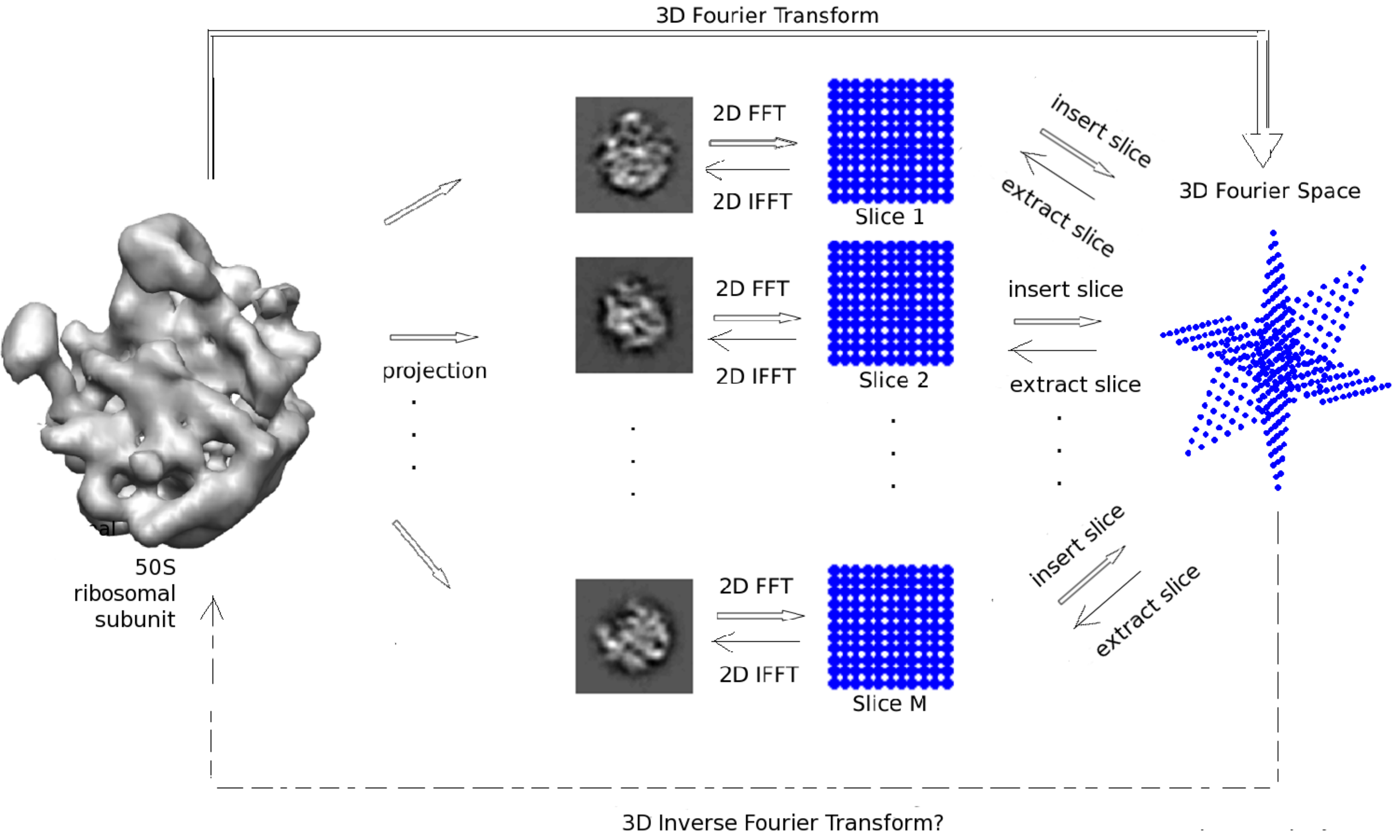}
\par\end{centering}

\caption{\label{fig:fourier_projection_slice_theorem}The Fourier projection-slice
theorem states that a slice extracted from the frequency
domain representation of a 3D volume yields the Fourier transform
of a projection of the volume in a direction perpendicular to the
slice. The volume
we show here is the 50S ribosomal subunit
used in our numerical experiments. The molecular surface was produced using the UCSF Chimera package \cite{Chimera} from the Resource for Biocomputing, Visualization, and Informatics at the University of California, San Francisco (supported by NIH P41 RR001081). }
\end{figure}

The Fourier projection-slice theorem plays a fundamental role in all 3D reconstruction algorithms independent of whether they are implemented in real space or in Fourier space \cite{vanheel_cryo-em}. The theorem states that a slice extracted from the frequency domain representation of a 3D map yields the 2D Fourier transform of a projection of the 3D map in a direction perpendicular to the slice (Figure \ref{fig:fourier_projection_slice_theorem}). It follows from the theorem that a reconstruction can be obtained by a 3D inverse Fourier transform from the Fourier domain which is filled in by the 2D Fourier slices. Although the continuous Fourier transform is a unitary linear transformation whose inverse equals its adjoint, the 3D discrete inverse Fourier transform of the slices does not equal its adjoint due to the non-uniform sampling in the frequency domain. Observe that the 3D Fourier space filled by 2D slices is denser at low frequencies and sparser at high frequencies. As a result, when the adjoint operator is applied to the slices, the low-frequency information of the macromolecule is overemphasized compared to the high frequency information, meaning that the inverse problem cannot be simply solved in this way. Instead, the solution to the linear inverse problem is either computed by applying a carefully designed weighted adjoint operator that addresses the non-uniform sampling \cite{GWBP_radermacher,GWBP_radermacher_weighted_1992,EWBP_1986,Penczek_GDFR}, or by using an iterative approach for inversion \cite{ART_Gordon1970471,ART_Marabini199853, SIRT_Gilbert1972}.

Many techniques have been developed to reconstruct a volume from images \cite{PawelA20101}. The Algebraic Reconstruction Technique (ART) \cite{ART_Gordon1970471,ART_Marabini199853} and the Simultaneous Iterative Reconstruction Technique (SIRT) \cite{SIRT_Gilbert1972} are algebraic approaches to find a 3D reconstruction such that its 2D re-projections are most similar to the input images in the least squares sense. The results of ART and SIRT are very accurate and they can incorporate additional constraints for the volume according to possible prior knowledge, such as positivity and smoothness. Another important advantage of ART and SIRT is that they are able to reconstruct from images with unevenly distributed viewing directions, which is usually the situation in cryo-EM since the macromolecules can assume various stable positions depending on their shape and adsorption properties \cite{vanheel_cryo-em,Frank1996}. However, ART and SIRT are extremely time-consuming if many iterations are needed for convergence. The filtered back-projection approaches, including the general weighted back-projection with exponent-based weighting function (WBP1) \cite{GWBP_radermacher,GWBP_radermacher_weighted_1992} and the exact filter weighted back-projection (WBP2) \cite{EWBP_1986} are considerably faster. However, the suitability of their weighting schemes depends on the distribution of the viewing directions, which can affect the precision of the reconstruction. The Gridding Direct Fourier Reconstruction (GDFR) \cite{Penczek_GDFR} is a relatively recent reconstruction technique. During the preprocessing stage, GDFR re-samples the 2D central slices onto 1D central radial lines to form a special structure of a non-uniform grid. Then the gridding weights are computed via a spherical Voronoi diagram. Finally, with the gridding weights, the numerical inverse Fourier transform is computed by the 3D gridding method. Although GDFR is both accurate and fast, it is limited to cases when there is no major gap among the viewing directions of the images since the proper gridding weights depend on the full coverage of Fourier space by the Fourier slices. Another direct Fourier inversion algorithm is the nearest neighbor  direct inversion reconstruction algorithm (4NN) \cite{4nn,ssn_4nn}. In the algorithm of 4NN, the 2D projections are first padded with zeros to four times the size, 2D Fourier transformed, and samples are accumulated within the target 3D Fourier volume using simple nearest neighbor interpolation. In the process, a 3D weighting function modeled on Bracewell’s “local density” \cite{4nn_weighting} is constructed and applied to individual voxels of 3D Fourier space to account for possible non-uniform distribution of samples. 4NN is even faster than GDFR and it is accurate when the sampling points are uneven in Fourier space. However, 4NN cannot avoid the projections whose Fourier transforms are close to gaps in Fourier space from receiving excessive weight.

\begin{figure}
\begin{centering}
\includegraphics[width=0.6\paperwidth]{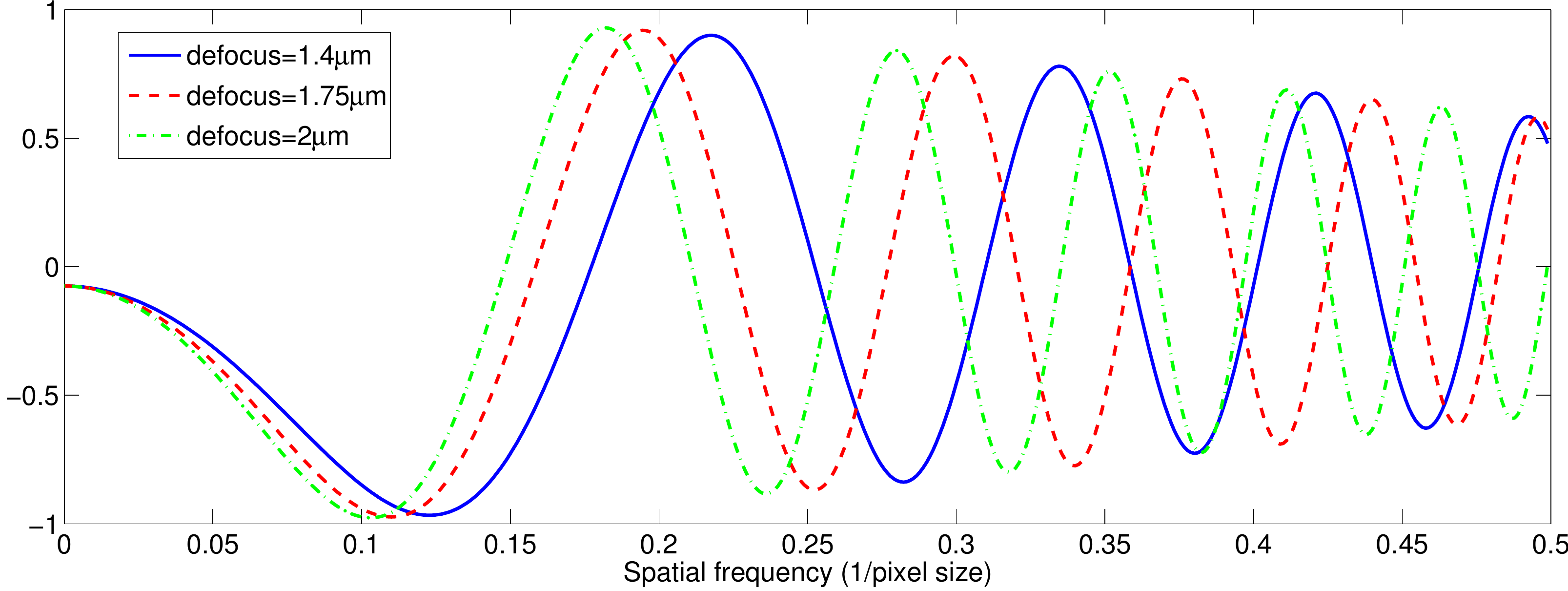}
\par\end{centering}

\caption{\label{fig:ctf}Three CTFs corresponding to different defocus values. The CTFs are generated
according to the defocus formula from page 81 of \cite{CTF_formula}. Parameters: $\alpha =.07$ (amplitude
contrast component), the electron wavelength $\lambda$= $2.51$ picometers, Cs=$2.0$ (spherical
aberration constant), B-factor=$100$, and defocus=$1.4\mu m$, $1.75\mu m$,
and $2\mu m$. Pixel size is $3.36$\AA.}
\end{figure}

Cryo-EM images, however, are not merely 2D mathematical projections of the macromolecule. During the imaging process, the objective lens of the electron microscope imposes a contrast transfer function (CTF) on a group of images \cite{Frank1996}. A CTF is approximated by a sinusoidal function in Fourier space depending on the magnitude of the frequency (Figure \ref{fig:ctf}), and it is also possible to improve the estimation of the CTF from the cryo-EM images themselves \cite{CTF, CTF_formula}. The CTF affects the acquired images through multiplication in the 2D Fourier domain, or equivalently, through a convolution in the real domain. The CTFs modulate the Fourier transform of true projections in a defocus-dependent way. A group of images taken using the same defocus setting is called a defocus group. One generally works at relatively large defocus values (up to $3\mu m$) to reduce the loss of low-frequency information of the images \cite{vanheel_cryo-em}. At large defocus values, the CTFs oscillate rapidly and decay exponentially in the high frequency domain (Figure \ref{fig:ctf}). The many zero crossings and fast decay of the CTFs cause the loss of information. Therefore, a good reconstruction must make use of images from different defocus groups, hoping that the information loss caused by the zero-crossings of one CTF would be filled by the information originating from images affected by other CTFs.

To reconstruct an undistorted volume from different defocus groups, one must carry out CTF correction of images or volumes. In the defocus groups approach, 3D CTF correction is applied to the reconstructed volume from each defocus group, then these volumes are combined to form a single, CTF-corrected volume \cite{Penczek_CTF1,CTF}. Another approach is to use CTF-corrected images or class averages for reconstruction. These two approaches take reconstruction and CTF-correction as two separate steps. In the step of reconstruction, one of the reconstruction algorithms (e.g. SIRT, WBP1, WBP2, GDFR, etc.) is used. In the step of CTF-correction, the Wiener filter is applied to find the least square solution to the problem of CTF-correction. 

In contrast with  these two approaches, it is possible to incorporate CTF corrections into the reconstruction algorithms in hope of a better merging of different defocus groups. In \cite{Penczek_CTF1}, Penczek et al. describe an algebraic method in real space to find a 3D reconstruction such that its 2D reprojections with CTFs are most similar to the input images in the least squares sense. Penczek et al. conclude that this approach outperforms the defocus groups approach. However, the algebraic method is time consuming. A similar idea is used in the refinement process by FREALIGN \cite{FREALIGN}, which creates a reconstruction by computing a least-squares fit to all the images with weights depending on the CTFs and the correlations between the images and the references. However, this reconstruction method is limited to the refinement process and cannot be applied for reconstructing an initial model. The direct Fourier inversion algorithm 4NN \cite{4nn,ssn_4nn} efficiently incorporates CTF corrections during the nearest neighbor interpolation using the Wiener filter methodology.

Here, we propose a fast and accurate Fourier-based Iterative Reconstruction Method (FIRM). FIRM is faster than ART and SIRT while maintaining their advantages over WBP1, WBP2, GDFR and 4NN: the results of FIRM are very accurate, it can incorporate prior knowledge, and does not require the viewing directions of the images to be evenly sampled. In addition, the time complexity of FIRM has the same order of magnitude as the fast algorithm 4NN and the actual running time differs only by a constant factor of around $5$.  Moreover, FIRM is also flexible to incorporate CTF corrections.

To derive the FIRM algorithm, we use a forward projecting model ${\bf b}={\bf A}({\bf V})+{\bf noise}$ according to the Fourier projection-slice theorem, where $\bf A$ is the forward projector, $\bf V$ is the unknown 3D density map of the macromolecule we are interested in, and ${\bf b}$ is the 2D Fourier transform of the noisy images. The reconstruction is obtained by minimizing the cost function $\rho\left({\bf V}\right)=\left\Vert {\bf b}-{\bf A}\left({\bf V}\right)\right\Vert ^{2}$. FIRM takes advantage of the Toeplitz structure of ${\bf A}^{*}{\bf A}$, where $\bf A^{*}$ is the adjoint operators of $\bf A$. The Toeplitz structure of the composition of the backward and forward projectors has already been successfully used for 2D reconstruction of an image from non-uniform Fourier-domain samples \cite{toep_2d,MRI,toeplitz,toeplitz_fessler}. We extend the usage of the Toeplitz structure to address the 3D reconstruction problem. Due to the Toeplitz structure of the operator ${\bf A^{*}A}$, it is equivalent to a convolution with a kernel. The kernel is precomputed using the non-uniform Fast Fourier Transform (NUFFT) \cite{dutt:1368,greengard:443,nufft_fessler} and is efficiently applied in an iterative process, such as the Conjugate Gradient (CG) method, to estimate the 3D map $\bf V$.

The outline of the paper is as follows. In section \ref{sec:introduction}, we introduce the necessary mathematical background concerning the reconstruction problem. In section \ref{sec:FIRM}, we demonstrate the key property that the composition of back-projection and forward-projection has a Toeplitz structure. We utilize this Toeplitz structure to accelerate the iterations in the CG method. Finally, numerical examples and concluding remarks are given in sections \ref{sec:numerical} and \ref{sec:conclusion}.

\section{Preliminaries}
\label{sec:introduction}

In this section, we provide the necessary mathematical background concerning the reconstruction
problem in cryo-EM and introduce notation used throughout this paper.

\subsection{Notation}
\label{sub:notation}

Scalars, indices and functions are denoted by non-boldface lowercase letters such as $x$, $i$, and $f$. Global constants are denoted by non-boldface uppercase letters such as $C$ and $N$. Boldface lowercase letters denote either vectors or
arrays, e.g., ${\bf a}=(a_{i_{1},i_{2},\ldots,i_{d}})$ with $i_{k}=1,2,\ldots,n_{k}$
and $k=1,2,\ldots,d$ is a $d$-dimensional array of size $n_{1}\times n_{2}\times\ldots\times n_{d}$. We refer to individual elements as either $a_{i_{1},i_{2},\ldots,i_{d}}$ or $a(i_1,i_2,\ldots,i_d)$. Operations such as
${\bf a}/2$ and ${\bf a}>0$ are considered component-wise. Matrices and operators are denoted by boldface uppercase letters, such as $\bf{A}$ and ${\bf F}$. The elements of a matrix ${\bf A}$ are denoted as either ${\bf A}(i,j)$ or $A_{i,j}$. The elements of a matrix ${\bf A}$ of multi-order $d>1$ are denoted as either $\bf A({\bf i}, {\bf j})$ or $A_{{\bf i}, {\bf j}}$, where  ${\bf i}$ and ${\bf j}$ are vector indices. The adjoint of a matrix
(or an operator) ${\bf A}$ is denoted as ${\bf A}^{*}$. The conjugate of a complex
number $z = a+\imath b$ is denoted as $\bar{z} = a-\imath b$. The absolute value of $z$
is denoted as $\left|z\right|=\sqrt{a^2+b^2}$. The inner product of two arbitrary $n_{1}\times n_{2}\times\ldots\times n_{d}$
arrays ${\bf a}$ and ${\bf b}$ is defined as $\left\langle {\bf a},{\bf b}\right\rangle =\sum_{i_{1},i_{2},\ldots,i_{d}}{a}_{i_{1},i_{2},\ldots,i_{d}}\bar{b}_{i_{1},i_{2},\ldots,i_{d}} = \sum_{{\bf i}} a_{{\bf i}} \bar{b}_{{\bf i}}$. We
omit the index and bounds of summation when these
are clear from the context.
The $\ell^2$ norm of $\bf{a}$ is denoted as $\left\Vert {\bf a} \right\Vert = \sqrt{\left \langle {\bf a}, {\bf a} \right \rangle}$.

\subsection{Fourier Transform Conventions}

The $d$-D Fourier transform ${\bf F}$ of a function $f:\mathbb{R}^{d}\rightarrow\mathbb{C}$ is
defined by
\[
({\bf F}f)(\boldsymbol{\omega})=\int_{\mathbb{R}^{d}}f({\bf x})\exp(-\imath\left \langle \boldsymbol{\omega},{\bf x} \right \rangle)\text{d}{\bf x},\quad \mbox{for } \boldsymbol{\omega}\in \mathbb{R}^d.
\]
Likewise, the $d$-D discrete Fourier Transform (DFT) ${\bf F}$ of an array ${\bf f}=(f_{{\bf n}})$,
where $-{N}/2\leq{\bf n}<{N}/2,$ is given by
\[
({\bf F} {\bf f})_{\bf k}=\sum_{-{N}/2\leq{\bf n}<{N}/2} f_{{\bf n}}\exp(-\imath2\pi\left \langle {\bf k},{\bf n}\right \rangle/N),\quad \mbox{for } -N/2 \leq {\bf k} < N/2.
\]

\subsection{Toeplitz Matrices and Circulant Matrices}
\label{sub:toeplitz_def}
We will show in section 3 that reconstructing a volume is equivalent
to solving a symmetric positive-definite Toeplitz system. In this subsection,
we introduce Toeplitz matrices and circulant matrices.
An $n\times n$ Toeplitz matrix is of the following form:
\[
{\bf T}_{n}=\left(\begin{array}{ccccc}
t_{0} & t_{-1} & \cdots & t_{2-n} & t_{1-n}\\
t_{1} & t_{0} & t_{-1} & \cdots & t_{2-n}\\
\vdots & t_{1} & t_{0} & \ddots & \vdots\\
t_{n-2} & \cdots & \ddots & \ddots & t_{-1}\\
t_{n-1} & t_{n-2} & \cdots & t_{1} & t_{0}
\end{array}\right),
\]
i.e., ${\bf T}_n(i,j)=t_{i-j}$ and ${\bf T}_{n}$ is constant along its diagonals.

A circulant matrix is a Toeplitz matrix of the form:
\[
{\bf C}_{n}=\left(\begin{array}{ccccc}
c_{0} & c_{n-1} & \cdots & c_{2} & c_{1}\\
c_{1} & c_{0} & c_{n-1} & \cdots & c_{2}\\
\vdots & c_{1} & c_{0} & \ddots & \vdots\\
c_{n-2} & \cdots & \ddots & \ddots & c_{n-1}\\
c_{n-1} & c_{n-2} & \cdots & c_{1} & c_{0}
\end{array}\right),
\]
i.e., ${\bf C}_n(i,j) = c_{i-j}$ where $c_{-k}=c_{n-k}$ for $1\leq k\leq n-1.$ Note that ${\bf C}_{n}$
is completely determined by its first column. It is well-known that circulant
matrices are diagonalized by the Fourier matrix ${\bf F}_{n}$ \cite{toeplitz_strange}, i.e.,
\begin{equation}
{\bf C}_{n}={\bf F}_{n}^{*}{\bf \Lambda}_{n}{\bf F}_{n},\label{eq:circ1}
\end{equation}
where ${\bf F}_{n}(j,k)=\frac{1}{\sqrt{n}}\exp\left(2\pi\imath jk/n\right),$
and ${\bf \Lambda}_{n}$ is a diagonal matrix. It follows immediately from (\ref{eq:circ1})
that the diagonal entries of ${\bf \Lambda}_{n}$, namely,
the eigenvalues of ${\bf C}_{n}$ can be obtained in $O\left(n\log n\right)$
operations using the Fast Fourier Transform (FFT) of the first column of ${\bf C}_{n}$. Once ${\bf \Lambda}_{n}$
is obtained, the matrix-vector product ${\bf C}_{n}{\bf y}$
can be computed efficiently by two FFTs in $O\left(n\log n\right)$
operations using (\ref{eq:circ1}) for any vector ${\bf y}$.

Similarly, we can define $\bf n$-by-$\bf n$ Toeplitz matrices and circulant matrices of multi-order $d$. The property (\ref{eq:circ1}) can be generalized to multi-order circulant matrices.

\section{A Fourier-based Approach for 3D Reconstruction}
\label{sec:FIRM}

\subsection{The Forward-Projector $\bf A$ }
In cryo-EM, the structure of a molecule is described by the molecule's
electric potential function $\varphi({\bf x}),$ where ${\bf x}=(x_{1},x_{2},x_{3})$
is in $\mathbb{R}^{3}.$ In a cryo-EM experiment, the macromolecules
are assumed to be identical with different orientations. We use $\bf R$ to denote
the rotation of each molecule, where $\bf R$ is an element of the rotation group ${\bf SO}(3)$.
The projection image of a molecule with orientation $\bf R$ is given
by
\[
({\bf P_{R}}\varphi)(x_{1},x_{2})=\int_{-\infty}^{\infty}\varphi_{\bf R}(x_{1},x_{2},x_{3})\text{d}x_{3},
\]
where $\varphi_{\bf R}({\bf x})=\varphi({\bf R}^{-1}{\bf x})$ is the electric
potential of the molecule after a rotation by $\bf R$. Note that for cryo-EM images, pose parameters include
both translations and rotations. Given the translations,
the images are re-shifted to their centers. Therefore, here we consider
the reconstruction problem for centered images given the rotational information.

With the above definitions of the Fourier transform and the projection,
one can verify the following theorem, known as the Fourier projection-slice theorem (page 11 in \cite{natterer_fourier_proj_slice}):
\begin{equation}
({\bf FP_{R}}\varphi)(\omega_{1},\omega_{2})=({\bf F}\varphi_{\bf R})(\omega_{1},\omega_{2},0).\label{eq:continuous_fps}
\end{equation}
The theorem states that the 2D Fourier transform of a projection of
an object $\varphi$ equals to one central slice of the 3D Fourier transform
of the object $\varphi$, where the projection is taken in a direction
perpendicular to the slice (Figure \ref{fig:fourier_projection_slice_theorem}).

It is important to realize that in practice the molecule's electric
potential function $\varphi$ is of limited spatial extent. On the other hand, numerically it is only possible to compute a finite discrete Fourier transform of $\varphi$. It is well known that a function with compact support cannot
have compactly supported Fourier transform unless it is identically
zero. However, this constraint is easily overcome for any finite accuracy \cite{signal_analysis}.
In this paper, the potential functions $\varphi$ are assumed to be essentially
band-limited to a ball and essentially space-limited to a cube. A
ball in the Fourier domain is a natural choice due to the radial symmetry of the CTFs and isotropic treatment
of orientations of cryo-EM images. We sample the continuous function
$\varphi$ on a Cartesian grid $\left\{ {\bf n}:\,{\bf n}\in\mathbb{Z}^{3},\,-{ N}/2\leq{\bf n}<{ N}/2\right\} $
to obtain a volume ${\bf V}\left({\bf n}\right)=\varphi\left({\bf n}a\right),$ where
$a\in\mathbb{R}_{+}$ is the grid spacing, and ${ N}\in\mathbb{Z}_{+}$
is large enough to cover the support of $\varphi$. According
to the sampling theorem, we further assume the Nyquist frequency $1/\left(2a\right)$
is no smaller than half the essential bandwidth of the function $\varphi$.  With
the above assumptions, the Fourier projection-slice theorem has the
following discretized version:
Given a volume $\bf V$ of size $N\times N\times N$ with the above assumptions,
and a projection's orientation ${\bf R}\in{\bf SO}(3)$, define the frequency
on the Cartesian grid of a central slice as $\boldsymbol{\omega}=(\omega_{1},\omega_{2})=2\pi(k_{1},k_{2})/N,$
where $k_{1},k_{2}\in\mathbb{Z}$. The Fourier projection-slice theorem (\ref{eq:continuous_fps}) implies that
the Fourier coefficient at $\boldsymbol{\omega}$ on the slice is approximated by
\begin{gather}
\sum_{{ -{N}/2\leq{\bf n}<{ N}/2}}{\bf V}_{{\bf n}}\exp\left(-\imath\cdot\left \langle{\bf n},{\bf R}^{-1}\left(\omega_{1},\omega_{2},0\right)\right \rangle\right).
\end{gather}
In particular if $\left\Vert \boldsymbol{\omega}\right\Vert >\pi$, then the Fourier coefficient at $\boldsymbol{\omega}$ is approximately zero.

Importantly, cryo-EM images are not true projections of a macromolecule
because of the effects of the CTFs \cite{Frank1996}. Mathematically, a CTF is
defined as a function in the Fourier domain, which can be approximated
by a sinusoidal function depending on the magnitude of the frequency. A CTF has the following form: 
\begin{equation*}
\text{CTF}\left(r\right)=\sin\left(-\pi\cdot\left(\text{defocus}\cdot r^{2}-\text{Cs}\cdot\lambda^{3}\cdot r^{4}/2\right)-\text{A}\right)\cdot\exp\left(-\left(\frac{r}{2\cdot\text{B factor}}\right)^{2}\right) ,
\end{equation*}
where $r$ is the magnitude of the frequency, Cs is the spherical aberration constant in mm, $\lambda$ is the electron wavelength in picometers, and A is amplitude contrast. A cryo-EM image $\bf I$ is the result of convolving the
true projection $\bf J$ with a point spread function, where the point spread function is the inverse Fourier transform of the CTF $h$. Thus, following the convolution
theorem, ${\bf F(I)}={\bf F(J)}h$, where $\bf F$ is the Fourier transform operation.

Denote a CTF as a function $h:\mathbb{R}_{+} \rightarrow \mathbb{R}$, then  according to (3), a Fourier slice affected by a CTF is approximated by
\begin{gather}
\sum_{-{ N}/2\leq{\bf n}<{ N}/2}{\bf V}_{{\bf n}}\exp\left(-\imath\cdot\left \langle{\bf n},{\bf R}^{-1}\left(\omega_{1},\omega_{2},0\right)\right \rangle\right)h(\left\Vert \boldsymbol{\omega}\right\Vert ).
\end{gather}

With the knowledge of the Fourier projection-slice theorem and the CTFs, it
is natural to define a forward-projector which projects a volume $\bf V$
to obtain Fourier slices modulated by CTFs.

Consider a volume $\bf V$ of size $N\times N\times N$, $M$ images with corresponding CTFs
$(h_{1},h_{2},\cdots,h_{M})$, and rotations ${\bf R}_{1},{\bf R}_{2},\ldots,{\bf R}_{M}\in{\bf SO}(3)$.
For each central slice (2D Fourier transform of images) consider the Cartesian coordinates $\boldsymbol{\omega}_{k_{1},k_{2}}=(\omega_{k_{1}},\omega_{k_{2}})=2\pi(k_{1},k_{2})/N,$
where $k_{1},k_{2}\in\mathbb{Z}$. We define a forward-projector $\bf A$ which
projects a volume ${\bf V}$ to obtain $M$ truncated Fourier slices corresponding
to the images as
\begin{equation}
\left({\bf A}\left({\bf V}\right)\right)\left(k_{1},k_{2},m\right)=\sum_{{\bf n}}{\bf V}_{{\bf n}}\exp\left(-\imath\cdot \left \langle {\bf n} , {\bf R}_{m}^{-1}\left(\omega_{k_{1}},\omega_{k_{2}},0\right) \right \rangle \right)\cdot h_{m}\left(\left\Vert \boldsymbol{\omega}_{k_{1},k_{2}}\right\Vert \right),\label{eq:forward-projector}
\end{equation}
where $m=1,\ldots,M$ is the index of an image, and ${\bf k}=(k_{1},k_{2})$ satisfies
the condition
\begin{equation}
\left\Vert {\bf k} \right\Vert\leq N/2 \text{      (inside a ball in the Fourier domain).}\label{eq:truncation}
\end{equation}
The condition (\ref{eq:truncation}) is based on the assumption that the function $\varphi$ corresponding to the volume $\bf V$
is essentially band-limited to a ball in the Fourier domain.
With this definition, the imaging process is modeled as
\begin{equation}
{\bf b}={\bf A(V)}+{\bf noise},\label{eq:measurement}
\end{equation}
where ${\bf b}$ is formed by the 2D discrete Fourier transform of the noisy images and restricting only to frequencies that satisfy (\ref{eq:truncation}). 

\subsection{The Back-projector ${\bf A}^{*}$ and the Toeplitz Structure of \textmd{${\bf A}^{*}{\bf A}$}}

The back-projector ${\bf A}^{*}$ is the adjoint operator of $\bf A$.
Note that ${\bf A}^{*}$ is not equivalent to the inverse of $\bf A$ because
of the non-uniform spacing of frequencies. Let ${\bf g}$ be an arbitrary collection of $M$ truncated slices, that is, ${\bf g} =(g_{k_{1},k_{2},m})$, with $1\leq m\leq M$, and ${\bf k}=(k_{1},k_{2})$ satisfy (\ref{eq:truncation}). From the definition of ${\bf A}^{*}$: $\left\langle {\bf A(V)},{\bf g}\right\rangle =\left\langle {\bf V},{\bf A}^{*}({\bf g})\right\rangle $,
we obtain
\begin{equation}
({\bf A}^{*}{\bf g})({\bf n})=\sum_{m=1}^{M}\sum_{{\bf k}}\exp\left(\imath\cdot\left\langle{\bf n},{\bf R}_{m}^{-1}\left(\omega_{k_{1}},\omega_{k_{2}},0\right)\right\rangle\right)\cdot h_{m}\left(\left\Vert \boldsymbol{\omega}_{k_{1},k_{2}}\right\Vert \right)\cdot g_{k_{1},k_{2},m}.\label{eq:backprojection}
\end{equation}
The operator ${\bf A}^{*}{\bf A}$ is then given by
\begin{eqnarray}
&&{\bf A}^{*}{\bf A}({\bf V})({\bf n})\nonumber\\
& = & \sum_{{\bf l}}{\bf V}_{{\bf l}}\sum_{m=1}^{M}\sum_{{\bf k}}\exp\left(\imath\cdot\left \langle{\bf n}-{\bf l},{\bf R}_{m}^{-1}\left(\omega_{k_{1}},\omega_{k_{2}},0\right)\right \rangle\right)\cdot h_{m}\left(\left\Vert \boldsymbol{\omega}_{k_{1},k_{2}}\right\Vert \right)^{2}\label{eq:toeplitz}\\
& = & (\bf Ker\star {\bf V})({\bf n}),\label{eq:toep2}
\end{eqnarray}
where the ``convolution kernel'' ${\bf Ker}$ is defined as
\begin{equation}
{\bf Ker}({\bf n})=\sum_{m=1}^{M}\sum_{{\bf k}}\exp\left(\imath\cdot\left \langle{\bf n},{\bf R}_{m}^{-1}\left(\omega_{k_{1}},\omega_{k_{2}},0\right)\right \rangle\right)\cdot h_{m}\left(\left\Vert \boldsymbol{\omega}_{k_{1},k_{2}}\right\Vert \right)^{2},\label{eq:kernel}
\end{equation}
and $-{ N}<{\bf n}=\left(n_{1},n_{2},n_{3}\right)<{ N}.$

From (\ref{eq:toeplitz}) we observe the Toeplitz structure of ${\bf A}^{*}{\bf A}$, which is a Toeplitz matrix of multi-order $3$ and
of size ${\bf N}$-by-${\bf N}$ (see section \ref{sub:toeplitz_def}). The $({\bf n},{\bf l})$ entry of ${\bf A}^{*}{\bf A}$ only depends on ${\bf n}-{\bf l}$, that is, $\left({\bf A}^{*}{\bf A}\right)\left({\bf n},{\bf l}\right) = \left({\bf A}^{*}{\bf A}\right)\left({\bf n}-{\bf l}\right).$
In fact, from (\ref{eq:toeplitz}), we have
\begin{equation}
\left({\bf A}^{*}{\bf A}\right)\left({\bf n},{\bf l}\right)=\sum_{m=1}^{M}\sum_{{\bf k}}\exp\left(\imath\cdot \left \langle \left({\bf n}-{\bf {\bf l}}\right), {\bf R}_{m}^{-1}\left(\omega_{k_{1}},\omega_{k_{2}},0\right)\right\rangle\right)\cdot h_{m}\left(\left\Vert \boldsymbol{\omega}_{k_{1},k_{2}}\right\Vert \right)^{2}.\label{eq:toep3}
\end{equation}
The Toeplitz structure allows us to rewrite (\ref{eq:toeplitz}) as (\ref{eq:toep2}), that is,
as a convolution of the volume $\bf V$ with the kernel $\bf Ker$, or simply summarized as
\begin{equation}
{\bf A}^{*}{\bf A}({\bf V})=\bf Ker\star {\bf V}.\label{eq:property}
\end{equation}

The circular convolution theorem tells us that the Fourier transform
of a convolution equals the product of the Fourier transforms.
Consider an $n$-by-$n$ Toeplitz matrix ${\bf T}_{n}$ and an arbitrary
$n$-vector ${\bf v}$. The matrix-vector multiplication ${\bf T}_{n}{\bf v}$
can be computed by 1D FFTs by first embedding ${\bf T}_n$ into a $2n$-by-$2n$
circulant matrix \cite{toeplitz_strange}, i.e.,
\[
\left[\begin{array}{cc}
{\bf T}_{n} & {\bf U}_{n}\\
{\bf L}_{n} & {\bf T}_{n}
\end{array}\right]\left[\begin{array}{c}
{\bf v}\\
{\bf 0}
\end{array}\right]=\left[\begin{array}{c}
{\bf T}_{n}{\bf v}\\
{\bf L}_{n}{\bf v}
\end{array}\right],
\]
where ${\bf U}_{n}$ and ${\bf L}_{n}$ are $n$-by-$n$ matrices designed in a way that ensures that the $2n$-by-$2n$ matrix is circulant.
Then, the multiplication is carried out by FFTs using the decomposition
(\ref{eq:circ1}) while ignoring the bottom half of the output vector (i.e., ${\bf L}_n {\bf v}$). The matrix-vector multiplication thus requires
$O\left(2n\,\log\left(2n\right)\right)$ operations. Similarly, due
to the Toeplitz structure of ${\bf A}^{*}{\bf A}$, the matrix-vector multiplication ${\bf A}^{*}{\bf A}\left({\bf V}\right)$ is a three-dimensional convolution that can
be computed using 3D FFTs by embedding the ${\bf N}$-by-${\bf N}$ matrix
${\bf A}^{*}{\bf A}$ into a $2{\bf N}$-by-$2{\bf N}$ circulant matrix ${\bf C^A}$
of multi-order $3$, and then carrying out the computation by using
the 3D version of (\ref{eq:circ1}) for decomposing ${\bf C^A}$ (see section \ref{sub:toeplitz_def}). Using the
property (\ref{eq:property}), it can be verified that the first {}``column''
of ${\bf C^A}$ is
\begin{equation}
C^{\bf A}_{{\bf i},{\bf 1}}={\bf Ker}\left(c\left(i_{1}\right),c\left(i_{2}\right),c\left(i_{3}\right)\right),\label{eq:embedding_circulant}
\end{equation}
where ${ 1}\leq{\bf i}=\left(i_{1},i_{2},i_{3}\right)\leq2{ N}$
and the function  $c$ is defined as
\begin{equation}
c\left(i\right)=\begin{cases}
i, & 1\leq i\leq N\\
1, & i=N+1\\
i-2N, & N+2\leq i\leq2N
\end{cases}.\label{eq:embedding_index}
\end{equation}
The computation of ${\bf A}^{*}{\bf A}\left(\bf V\right)$ thus requires $O\left(8N^{3}\,\log\left(8N^{3}\right)\right)$ operations.

\subsection{The Conjugate Gradient (CG) Method}
\label{subsec:cg}
We reconstruct the volume by minimizing the cost function
\begin{equation}
\rho\left({\bf V}\right)=\left\Vert {\bf b}-{\bf A}\left({\bf V}\right)\right\Vert ^{2},\label{eq:cost1}
\end{equation}
where ${\bf b}$ includes all measured Fourier slices, $\bf A$ is
the forward-projector, given by (\ref{eq:forward-projector}), and $\bf V$
is the unknown volume. A solution to the minimization problem can
be found by setting to zero the derivative of $\rho$ with respect to $\bf V$, yielding
\[
{\bf A}^{*}{\bf A}\left({\bf V}\right)={\bf A}^{*}{\bf b}.
\]
Since ${\bf A}^{*}{\bf A}$ is symmetric positive-semidefinite,
we can apply CG to find the minimizer to the cost function (\ref{eq:cost1}).

To be used in an iterative method, the operation of $ {\bf A}^{*}{\bf A}$
must be extremely efficient. However, applying $\bf A$ and ${\bf A}^{*}$ separately
at each iteration is time-consuming since the time cost of one application
of either $\bf A$ or ${\bf A}^{*}$ is equivalent to the cost of one application of
NUFFT, whose time complexity is $O\left(MN^{2}\max\left(\log M,\,\log\left(N^{2}\right)\right)\right)$  \cite{dutt:1368,greengard:443,nufft_fessler}. However, this efficiency problem can be overcome using
the property (\ref{eq:property}). 
Thus, applying ${\bf A}^{*}{\bf A}$ (or equivalently ${\bf Ker}$) to a vector
requires $O\left(8N^{3}\,\log\left(8N^{3}\right)\right)$ operations.
Note that both
${\bf Ker}$ and ${\bf A}^{*}{\bf b}$ can be precomputed using NUFFT only once before applying
the CG method. More details about the time complexity of the computation of ${\bf Ker}$ and ${\bf A}^{*}{\bf b}$ are provided in section \ref{sub:FIRM}.

\begin{figure}
\begin{centering}
\includegraphics[width=0.6\paperwidth]{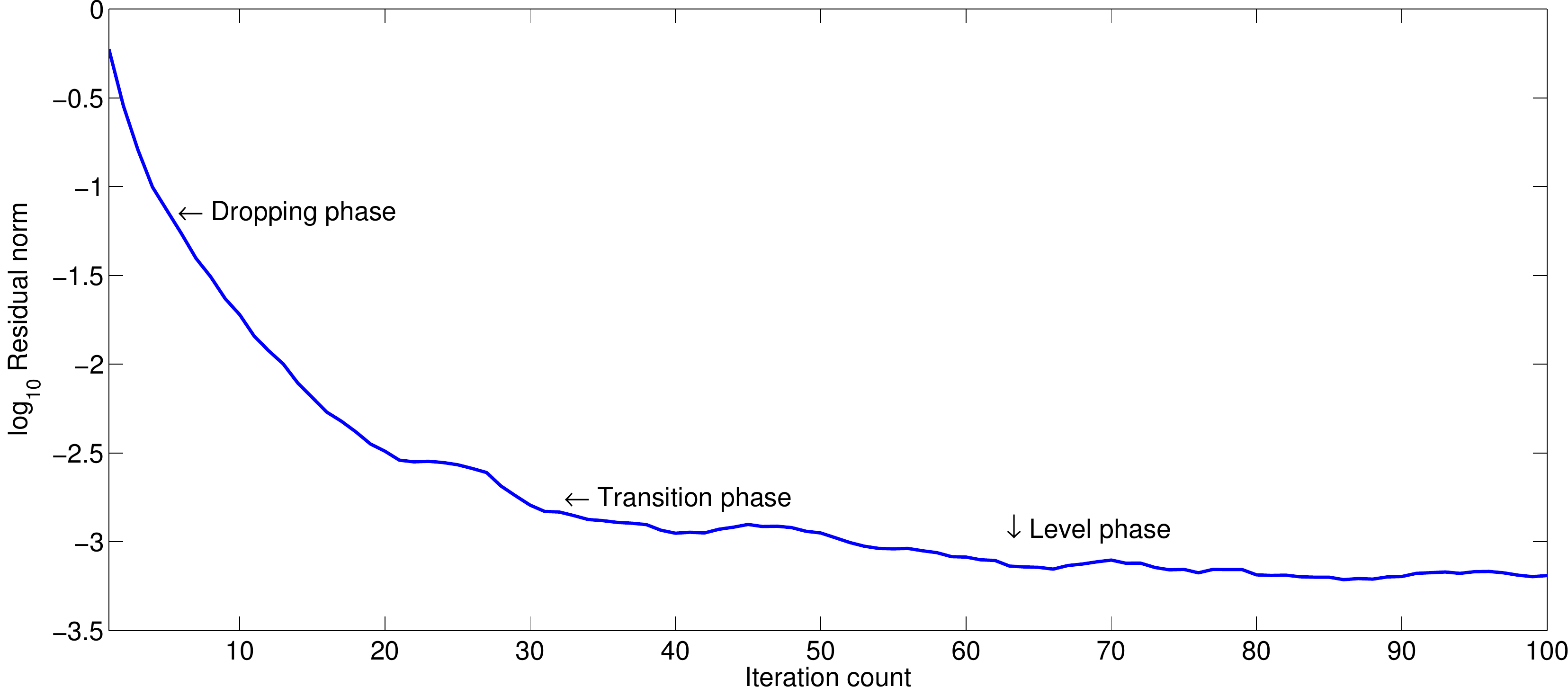}
\par\end{centering}

\caption{\label{fig:phase}Plot of residual norm (in $\log_{10}$  scale) versus iteration count for a typical CG iterative 3D reconstruction from images. The L-curve is divided into 3 segments, which correspond to “dropping phase”, “transition phase”, and “level phase” respectively.}
\end{figure}

The convergence rate of the CG method has been well studied (see \cite{Nocedal99}), and it depends on  how clustered the spectrum of ${\bf A}^{*}{\bf A}$ is. The output of the projector ${\bf A}$ has no high frequency information outside a ball in the Fourier domain, resulting in an extremely large condition number of ${\bf A}^{*}{\bf A}$ and very small eigenvalues in the spectrum. The ill-conditioning of the reconstruction problem causes the semi-convergence behavior \cite{semi_convergence_1,semi_convergence_2,semi_convergence_3}, which can be characterized as initial convergence toward the exact solution and later divergence. From the perspective of regularization, the updated estimated volume in each iteration is a regularized solution and the number of iterations plays the role of the regularization parameter. The iteration count controls the compromise between the signal-to-noise ratio (SNR) and the residual aliasing artifact. The basic principle of regularization is to “smooth” the solution by truncating or damping the small eigenvalue components.  The initial iterations pick up the eigenvalue components corresponding to the largest eigenvalues. As the iteration number increases, more and more small eigenvalues are captured and the degree of regularization decreases. As a result, the residual norm declines sharply at early stages of the iterative process before it levels off. For a well behaved reconstruction, the plot of the residual norm in $\log_{10}$  scale versus the iteration count generally exhibits an L-curve characteristic, as shown in Figure \ref{fig:phase}. In this manner, the iteration procedure can be divided into 3 phases. On the left of the L-curve, the residual norm declines very fast and we refer to it as a “dropping phase”; on the right side, the residual norm levels off and it is a “level phase”. The L-corner then represents the “transition phase”, where the noise and artifacts are usually well compromised. However, it is practically difficult to locate the point where SNR and artifacts are “optimally” compromised and “best” reconstructed volume quality is achieved. It is recommended to pick a number of reconstructed volumes in the transition phase for better presentation of the reconstruction results.

\subsection{Fourier-based Iterative Reconstruction Method (FIRM)}
\label{sub:FIRM}

The input to our reconstruction algorithm is the following:
\begin{enumerate}
\item $M$ projection images of size $N\times N$ denoted ${\bf I}_{1},{\bf I}_2,\ldots,{\bf I}_{M}$
of an unknown volume $\bf V$ of size $N\times N\times N$.
\item The orientations of the images ${\bf R}_{1},{\bf R}_{2},\ldots,{\bf R}_{M}\in{\bf SO}(3)$.
\item The CTFs $h_1, h_2, \ldots, h_M$.
\end{enumerate}
The description of FIRM is as follows:
\begin{enumerate}
\item Compute the 2D DFT of all images using 2D FFT. Truncate
the square Fourier slices and form the vector ${\bf b}$ in (\ref{eq:measurement}).
$\left(O\left(MN^{2}\log\left(N^{2}\right)\right)\right)$
\item Compute the back-projection ${\bf A}^{*}{\bf b}$ according to (\ref{eq:backprojection}) using NUFFT.
\newline $\left(O\left(MN^{2}\max\left(\log M,\,\log\left(N^{2}\right)\right)\right)\right)$
\item Compute the convolution kernel ${\bf Ker}$ according to (\ref{eq:kernel})
using NUFFT. 
\newline$\left(O\left(4MN^{2}\max\left(\log M,\,\log\left(N^{2}\right)\right)\right)\right)$
\item Use CG with input $\bf A^{*}{\bf b},\, {\bf Ker}$ and an initial guess (all-zero volume). The output is the estimated volume.
($O\left(8N^{3}\log\left(8N^{3}\right)\right)$ operations for each iteration)
\end{enumerate}
The running time of the algorithm is dominated by Steps 2-3. The time complexity of Step 2 for the NUFFT algorithm is obtained from \cite{dutt:1368,greengard:443,nufft_fessler}. In Step 3, although ${\bf Ker}$
is about $8$ times as large as the original volume (i.e., a factor of 2 in each dimension), by noting that ${\bf Ker}\left(-{\bf n}\right)=\overline{{\bf Ker}\left({\bf n}\right)}$, the time cost of computing ${\bf Ker}$ is about $4$ times the cost of computing the back-projection ${\bf A}^{*}{\bf b}$. Using
property (\ref{eq:property}), the computation in each iteration of CG in Step 4 is efficiently reduced to FFTs and matrix-matrix point-wise multiplication. Moreover, it is easy to parallelize the computation of $\bf A^*b$ and $\bf Ker$ by noting that both of them are summation over back-projection of a single projection. A MATLAB package for FIRM is available to download through the website \url{http://www.math.princeton.edu/~lanhuiw/software.html}.

\section{Numerical Results}
\label{sec:numerical}

We implemented FIRM
using the MATLAB programming Language. The NUFFT package provided by \cite{nufft_fessler}
 is used for precomputation of back-projections and convolution kernels. We compare FIRM with 4NN implemented within the framework of the SPARX image processing
system \cite{sparx}. The numerical experiment is performed on a machine with 2 Intel(R)
Xeon(R) CPUs X5570, each with 4 cores, running at 2.93 GHz. Both MATLAB and SPARX are limited to single core computations.

In the numerical experiment, a 50S ribosomal subunit volume of size $90\times 90\times 90$ (Figure \ref{fig:fourier_projection_slice_theorem}
Left) is used to generate projections. We use SPARX to generate a random conical tilt series consisting of 10,000 simulated projections. The tilt angle is fixed to be $60^\circ$ and the azimuthal angles are randomly sampled from the uniform distribution over $[0^\circ, 360^\circ]$. Thus there is a missing cone in the coverage of Fourier space by the slices. The purpose to use a random conical tilt series is that not only the accuracy of the reconstructed volumes excluding the missing cone can be studied, but also the artifacts of the volumes inside the missing cone can be observed. For a real dataset in cryo-EM, the geometry of the collected images cannot be totally controlled. Thus the artifacts of the reconstructions due to the uneven sampling are of interest \cite{uneven}.

The 10,000 projections are divided randomly to 3 defocus groups and filtered by the CTFs which are generated with parameters detailed in the caption of Figure \ref{fig:ctf}. We refer the CTF filtered projections as clean images. The noisy images are generated by  adding white Gaussian noise to the clean images. In this experiment, the SNR of the noisy image is set to 1. 

To evaluate the accuracy of the reconstructions, we use  the 3D Fourier
Shell Correlation (FSC) \cite{FSC}. FSC measures the normalized cross-correlation
coefficient between two 3D volumes over corresponding spherical shells
in Fourier space, i.e.,
\[
\text{FSC}\left(i\right)=\frac{\sum_{{\bf j} \in Shell_i}{\bf F}\left({\bf V}_{1}\right)\left({\bf j}\right)\cdot \overline{{\bf F}\left({\bf V}_{2}\right)\left({\bf j}\right)}}{\sqrt{\sum_{{\bf j} \in Shell_i}\left|{\bf F}\left({\bf V}_{1}\right)\left({\bf j}\right)\right|^{2}\cdot\sum_{{\bf j} \in Shell_i}\left|{\bf F}\left({\bf V}_{2}\right)\left({\bf j}\right)\right|^{2}}},
\]
where ${\bf F}\left({\bf V}_{1}\right)$ and ${\bf F}\left({\bf V}_{2}\right)$ are the Fourier transforms of volume ${\bf V}_1$
and volume ${\bf V}_2$ respectively,  the spatial frequency $i$ ranges from $1$
to $N/2-1$ times the unit frequency $1/(N\cdot \text{pixel size})$, and $ Shell_i:=\{{\bf j}:0.5+(i-1)+\epsilon\leq\left\Vert {\bf j}\right\Vert < 0.5 + i +\epsilon\}$ where $\epsilon =$1e-4. In this form, the FSC takes two 3D volumes
and converts them into a 1D array. For each reconstructed volume,
we measure its FSC against the clean 50S ribosomal subunit volume, that is, in our measurement ${\bf V}_1$ is the reconstructed volume, and ${\bf V}_2$ is the ``ground truth'' volume. In this case, FSC is also called Fourier Cross-Resolution (FCR) \cite{FCR}. To measure the accuracy of the reconstructed volumes outside and inside the missing cone respectively, we use modified FCR for the target Fourier volumes excluding or within the missing cone region.

\begin{figure}
\begin{centering}
\includegraphics[width=0.6\paperwidth]{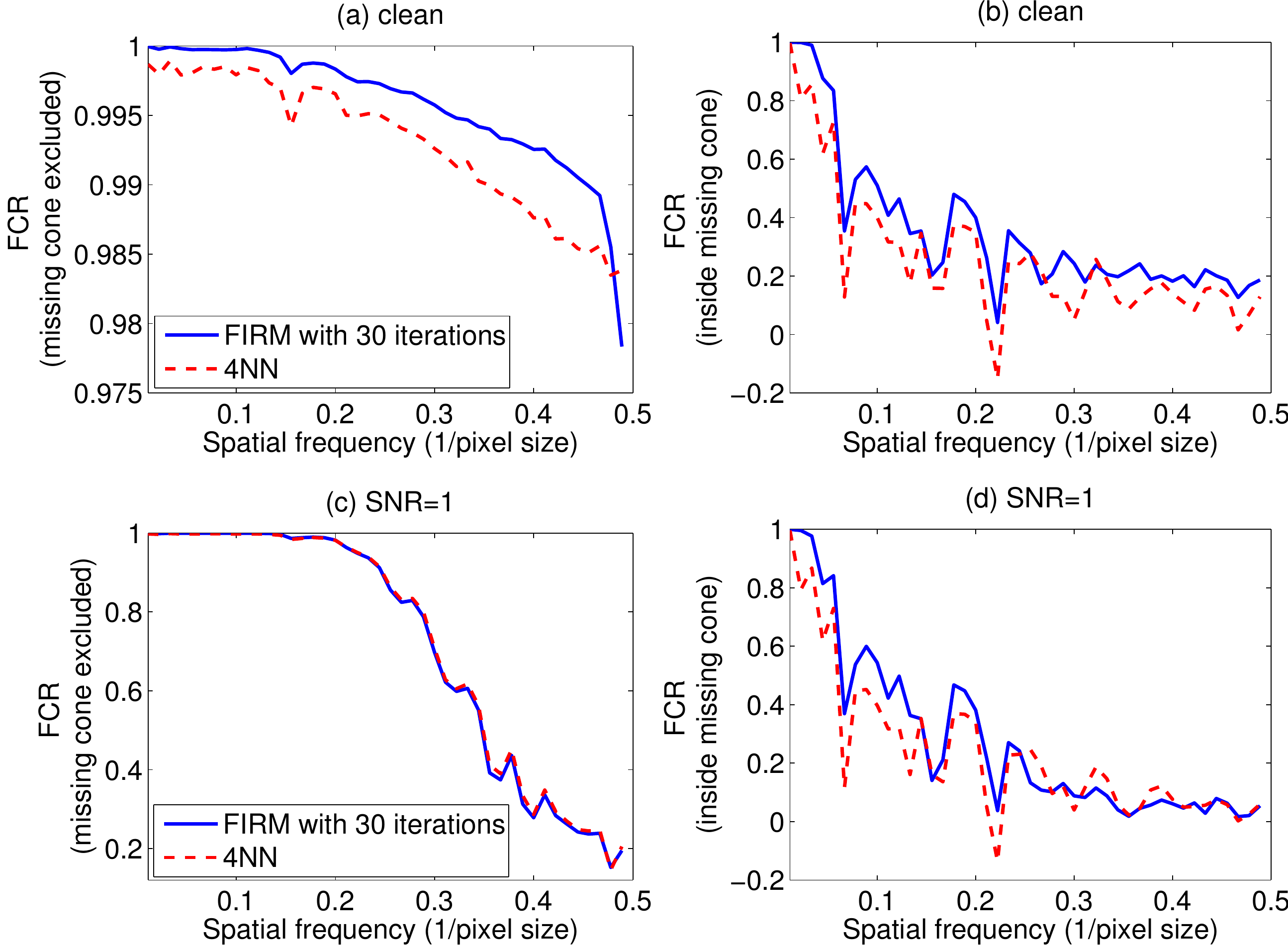}
\par\end{centering}

\caption{\label{fig:fsc}FCR of the reconstructed volumes excluding or within the missing cone region. (a) and (b) are comparison between reconstructions from clean images. (c) and (d) are comparison between reconstructions from noisy images with SNR=1. }
\end{figure}

FIRM is compared with other reconstruction algorithms (GDFR, SIRT, 4NN, etc). However we only report the comparison with 4NN since it performed best in terms of accuracy and running time  \cite{4nn,ssn_4nn}. The reconstructed volumes by FIRM are the estimations in $30$th iteration for both the clean and noisy image datasets, where the ``transition phase'' is in the L-curve (see details in section \ref{subsec:cg}). FIRM spent 4 seconds on FFTs of the images, 293 seconds on back-projection, 1143 seconds on computing the kernel, and 1 second on each CG iteration. The total time cost by FIRM is around 1470 seconds, which is about five times the time cost by 4NN (290 seconds). From Figure \ref{fig:fsc} we observe that the accuracy of the reconstructions by the two algorithms are almost the same excluding the missing cone. However, the measurement inside the missing cone demonstrates that there is less artifacts created by FIRM than by 4NN. 

\section{Summary and Discussion}
\label{sec:conclusion}
In this paper, we propose a fast and accurate Fourier-based iterative reconstruction method (FIRM) to reconstruct molecular structures from cryo-EM images. To study the imaging process in cryo-EM, we define a forward-projector $\bf A$ which converts a given volume to Fourier slices affected by CTFs. Therefore, the imaging model is $\bf{b}={\bf A(V)}+\bf{noise}$ where $\bf A$ is the forward-projector, $\bf V$ is the unknown volume we are interested in, and $\bf{b}$ is the measurement of the Fourier slices. To reconstruct the volume $\bf V$ from the measurement $\bf{b}$, CG is applied to estimate the reconstructed volume by minimizing $\rho\left({\bf V}\right)=\left\Vert {\bf b}-{\bf A}\left({\bf V}\right)\right\Vert ^{2}$. The solution is found by setting the derivative of $\rho$ to zero, yielding the equation
$
{\bf A}^{*}{\bf A}\left({\bf V}\right)={\bf A}^{*}{\bf b}.
$
$\bf A^*A$ has Toeplitz structure and thus ${\bf A}^{*}{\bf A}\left({\bf V}\right)={\bf Ker}\star {\bf V}$, where ${\bf Ker}$ is a convolution kernel. Using this property, which is key to our method, the computation of ${\bf A}^{*}{\bf A}\left({\bf V}\right)$ is reduced to FFTs and matrix-matrix point-wise multiplication according to the convolution theorem. As a result, the computation of each CG iteration is fast.

The main advantage of iterative methods (ART, SIRT and FIRM) in general is their applicability to diverse data collection geometries and to data with uneven distribution of projection directions. However, for ART and SIRT, the computational requirements are dominated by the back-projection steps and thus their running time exceeds that of other algorithms (WBP1, WBP2, GDFR and 4NN) for the typical number of required iterations (typically 10-200). Instead of back-projecting in each iteration, FIRM computes the back-projection only once in the pre-computation stage that also includes the computation of the kernel using NUFFT.  The most time cost by FIRM is thus at the stage of the preparation before CG. The time cost of iterations in FIRM is negligible compared to the pre-computation. 

The numerical experiments demonstrate that compared with 4NN, FIRM is fast and accurate, and it performs well in merging information from different defocus groups. Moreover, FIRM creates a satisfactory reconstruction in the case of a missing cone in Fourier space with less artifacts than 4NN. As of the running time, FIRM costs about five times the time used by 4NN, which is the  fast direct inversion algorithm in SPARX.
We remark that FIRM is flexible in the sense that it can incorporate other regularization terms that can be naturally formulated in Fourier space (e.g., damping high frequencies), or other prior knowledge about the volume, such as positivity constraints. In severe situations when images correspond to only a small number of views or when there are some gaps in Fourier space, regularization can alleviate the ill-conditioning of the problem. In \cite{WDFR} we demonstrated how the computational framework of FIRM can be modified to treat a regularization term that consists of the $\ell^1$ norm of the wavelet expansion coefficients of the volume. The purpose of such regularization terms is to promote sparsity in the expansion of the volume in the wavelet basis. We remark that other regularization terms involving the total variation functional \cite{TV} or tight frame expansions are also possible.
\section{Acknowledgements}
The authors would like to thank Fred Sigworth and C\'edric Vonesch for useful discussions and suggestions. The work was partially supported by Award Number R01GM090200 from the NIGMS. L. Wang was partially supported by Award Number DMS-0914892 from the NSF, Y. Shkolnisky was supported in part by Israel Science Foundation grant 485/10, and A. Singer was partially supported by the Alfred P. Sloan Foundation.

\bibliographystyle{plain}
\bibliography{cryo-EM_lanhui}

\end{document}